\documentclass{jgcc}

\pdfoutput=1

\usepackage{lastpage}

\jgccdoi{16}{1}{2}{12282}

\jgccheading{}{\pageref{LastPage}}{}{}{Sept.~16,~2023}{May.~15,~2024}{}

\usepackage{latexsym}
\usepackage{amsfonts}
\usepackage{hyperref}
\usepackage{xcolor}

\newtheorem{theorem}{Theorem}

\newtheorem{proposition}{Proposition}

\begin{document}

\author[]{Vladimir Shpilrain}
\address{Department of Mathematics, The City College of New York, New York,
NY 10031} \email{shpilrain@yahoo.com}

\title{On isomorphisms to a free group and beyond}

\begin{abstract}
The isomorphism problem for infinite finitely presented groups is probably the hardest among standard algorithmic problems in group theory. Classes of groups where  it has been completely solved are nilpotent groups, hyperbolic groups, and limit groups. In this short paper, we address the problem of isomorphism to particular groups, including free groups. We also address the algorithmic problem of embedding a finitely presented group in a given limit group.

\end{abstract}

\maketitle

\hfill{\small \it In memory of Ben Fine}

\section{Introduction}

The isomorphism problem has been completely solved in the class of finitely generated  nilpotent groups in \cite{Grune}.

Later, it was solved in the class of hyperbolic groups \cite{Sela} (torsion-free case), \cite{Guira} (general case), although it is difficult (if at all possible) to ``computerize" these algorithms, i.e., to code them in one of the known programming languages.

Then, the isomorphism problem was also solved in the class of limit groups (a.k.a. fully residually free groups) \cite{BumaginKM}.

In the class of finitely generated one-relator groups, although the isomorphism problem is still open in general, it has been settled for ``most" one-relator groups (in a precise formal sense) in \cite{KSS}. More specifically, for any $r \ge 2$, there is a subset $\mathcal{G}$ of elements of the free group $F_r$ such that: (1) $\mathcal{G}$ has asymptotic density 1 in $F_r$; (2) it is algorithmically possible to find out whether or not a given element $u \in F_r$ is in $\mathcal{G}$; (3) for any two elements $u, v \in \mathcal{G}$, it is algorithmically possible to find out whether or not two one-relator groups (with the relators $u$ and $v$, respectively) are isomorphic.

The ``next in line" class of groups where the isomorphism problem may be solvable is the class of finitely presented metabelian groups (see \cite{problems}, Problem (M1)), where ``most" algorithmic problems have solutions by now \cite{BCR}.

We note that the isomorphism problem has a reasonable chance to be solvable only in classes of groups where all groups have solvable word problem. This rules out, for example, the class of finitely presented solvable groups of derived length $\ge 3$ since this class has groups with unsolvable word problem \cite{OK}.

In this paper, we address an apparently easier problem of isomorphism to a particular group.
Using a simple trick, we establish here the following result that appears to be useful in some situations.

\begin{proposition}\label{main}
Let $G$ be a group with $n$ given generators. Suppose that $G$ has solvable word problem. Let $H$ be a finitely presented group, and suppose either $G$ or $H$ is Hopfian. If one can decide whether or not there is an epimorphism from $G$ onto $H$ and find it as an explicit map on the generators in case it exists, then one can decide whether or not $G$ is isomorphic to $H$.

\end{proposition}

Recall that a group is {\it Hopfian} if any {\it onto} endomorphism of this group is also one-to-one, i.e., is an automorphism. Note that in Proposition \ref{main} we do not require that $H$ has solvable word problem or that $G$ is finitely presented.

%Since finitely generated free solvable groups are Hopfian (in fact, residually finite \cite[41.52]{Neumann}) and have solvable word problem (in fact, solvable very efficiently, see \cite{Vershik}, \cite{Ushakov}), we have
%
%\begin{corollary}\label{solvable}
%Let $S_{n,d}$ be the free solvable group of rank $n$ and  derived length $d$. Then, given any recursively presented solvable group $H$ of derived length $d$ and given $n$ generators of $H$, it is algorithmically possible to find out whether or not $H$ is isomorphic to $S_{n,d}$.
%
%\end{corollary}
%
%This corollary actually holds for any free polynilpotent group, but there are probably very few people left who remember what this is.
%
%We note again that we do not {\it a priori} require that $H$ has solvable word problem. On the other hand, the condition on $H$ to have the same number $n$ of generators as the group $S_{n,d}$ does, is very restrictive.

Our main goal actually was to address the problem of isomorphism to the (absolutely) free group $F_n$ of rank $n$. There is a classical result of Adyan \cite{Adyan} saying that given an arbitrary (finitely presented) group $B$, there is no algorithm that would decide, given any (finitely presented) group $G$, whether or not $G$ is isomorphic to $B$.
However, if we require solvability of the word problem in $G$, then the problem of isomorphism of $G$ to the free group $F_n$ becomes algorithmically solvable:

\begin{theorem}\label{free}
Let $G$ be a finitely presented group with a given algorithm for solving the word problem in $G$. Then, for any given $n \ge 1$, it is algorithmically possible to find out whether or not $G$ is isomorphic to a free group of rank $n$.

\end{theorem}

%To the best of our knowledge, this result has not been explicitly mentioned anywhere in the literature, although
There is a ``detour" that leads to this result, see \cite[Corollary 4.3]{GrovesW} for an explicit mention of this result. Specifically, there is an algorithm that, given a finitely presented group $G$ with solvable word problem, decides whether or not $G$ is a limit group \cite{Groves}. If not, then $G$ cannot be isomorphic to a free group because any finitely generated free group is a limit group. If $G$ is a limit group, then one can use an algorithm, due to \cite{BumaginKM}, that decides if there is an isomorphism between two limit groups.

Our proof is more straightforward, but it still uses a ``big gun", namely Razborov's work on solving (systems of) equations in a free group.

It appears that solvability of equations in groups should inevitably be an important ingredient in any solution of the isomorphism problem for infinite groups. However, this is typically not enough. In our proof of Theorem \ref{free}, we actually establish an isomorphism (or non-isomorphism) of the group $G$ to a subgroup of a given fixed finitely generated free group, and then we use the fact that every nontrivial subgroup of a free group is itself free. This is not the case with hyperbolic groups, say; moreover, a finitely generated subgroup of a hyperbolic group may not even be finitely presented, and this makes our method inapplicable in that situation. One class of groups where our method does work is the class of limit groups since every finitely generated subgroup of a limit group is a finitely presented limit group. Also, finitely generated limit groups are Hopfian because they are residually free and therefore residually finite.  The following result may be of interest:

\begin{theorem}\label{limit}
Let $G$ be a finitely presented group with a given algorithm for solving the word problem in $G$. Let $H$ be a limit group. Then it is algorithmically possible to find out whether or not $G$ can be embedded in $H$.
\end{theorem}

\section{Proof of Proposition \ref{main}}

Let $g_1,\ldots, g_n$ be the given generators of the group $G$, and $h_1,\ldots, h_n$ generators of the group $H$. Needless to say, if there is no epimorphism from $G$ onto $H$, then $G$ and $H$ are not isomorphic.

Now suppose the map $\varphi : g_i \to h_i$ can be extended to an epimorphism from $G$ onto $H$. Then run two algorithms in parallel:

\medskip

\noindent {\bf 1.} Algorithm $\mathcal{A}$ will detect non-isomorphism by looking for an element in the kernel of $\varphi$. To that effect, it goes over nontrivial elements of $G$ one at a time (this is possible since the word problem in $G$ is solvable) and  checks if $\varphi$ takes them to the trivial element of $H$.

Here the reader may say: wait, you do not require that the word problem in $H$ is solvable. Indeed, but here we only need the ``yes" part of the word problem (i.e., detecting that the element is trivial), and this part works in any recursively presented group. Specifically, to detect that $w=1$ one can go over all finite products of conjugates of defining relators and (graphically) compare them to $w$.

We note that if the kernel of $\varphi$ is nontrivial, then $H$ is isomorphic to a proper factor group of $G$ and therefore cannot be isomorphic to $G$ since we assumed that either $G$ or $H$ was Hopfian.

\medskip

\noindent {\bf 2.} Algorithm $\mathcal{B}$ will detect isomorphism by looking for a map $\psi$, given on the generators $h_i$ of $H$, such that $\psi(\varphi(g_i))=g_i$ for all generators $g_i$ of the group $G$. To that effect, $\mathcal{B}$ will go over $n$-tuples $(y_1, \ldots, y_n)$ of elements of $G$, one at a time, and define $\psi$ by $\psi(h_i)=y_i$.

First check if $\psi$ is a homomorphism by computing $\psi(r_j)$ for every defining relator $r_j$ of the group $H$ and checking if $\psi(r_j)=1$. This is possible since $G$ has solvable word problem, although we do not really need this because again, here we only need the ``yes" part of the word problem.

If $\psi$ is a homomorphism, then just check if $\psi(\varphi(g_i))=g_i$ for all $g_i$, again using the ``yes" part of the word problem in $G$. If $H$ is isomorphic to $G$, then eventually a map $\psi$ like that will be found.

\medskip

Eventually one of the algorithms, $\mathcal{A}$ or $\mathcal{B}$, will stop and give an answer. $\Box$

We note that the only place in the proof where we used solvability of the word problem in $G$ was where we were trying to detect non-isomorphism by looking for a nontrivial element in the kernel of $\varphi$.

%Corollary \ref{solvable} does not follow directly from Proposition \ref{main} because in Proposition \ref{main}, one of the groups has to be finitely presented, which is not the case in Corollary \ref{solvable}. (We note, in passing, that the group $S_{n,d}$ is not finitely presented for any $n \ge 2, ~d \ge 2$.)
%
%On the other hand, in Corollary \ref{solvable} we require that the groups $S_{n,d}$ and $H$ have the same number $n$ of generators, and this simplifies the situation greatly.

\section{Proof of Theorem \ref{free}}

Let $G$ have $m$ generators.
First we note that if $n>m$, then $G$ cannot be isomorphic to a free group of rank $n$, so we assume that $n \le m$ from now on.

Let $g_1,\ldots, g_m$ be the given generators of the group $G$, and let $r_1,\ldots, r_s$ be all defining relators of $G$. Let $F_n$ be a free group of rank $n$, and let $\alpha : g_i \to x_i$ for some $x_i \in F_n$, $i=1, \ldots, m$. This map extends to a homomorphism $\alpha : G \to F_n$ if and only if $\alpha(r_j)=1$ for all $j=1, \ldots, s$. This translates into a system of $s$ equations in the group $F_n$.

%The set of all solutions of such a system of equations is recursively enumerable. Denote this set of $m$-tuples of elements of $F_n$ by $S$.
First, we will run Razborov's algorithm $\mathcal{R}$ \cite{Razborov} to see if this system of equations has a
solution tuple $(a_1,\ldots, a_m)$ that generates a free subgroup of rank $r \ge n$ in $F_n$; in other words, if there is an epimorphism of $G$ onto a free group of rank $r \ge n$. Denote this free group by $H_r$ (recall that every nontrivial subgroup of a free group is free). If the system has no such solutions, then $G$ is not isomorphic to a free group of rank $n$.

If there is an epimorphism of $G$ onto a free group $H_r$, then there is also an epimorphism of $G$ onto a free subgroup $K_n \le H_r$ of rank $n \le r$. To find an explicit epimorphism of $G$ onto $K_n$ (as a map on the generators), one can first select generators of $K_n$ and then find an explicit epimorphism of $H_r$ onto $K_n$ by using, say, Nielsen's method, see e.g. \cite{Lyndonbook}.

After one finds an epimorphism of $G$ onto $K_n$, Proposition \ref{main} applies (since any finitely generated free group is Hopfian), and this completes the proof. $\Box$

We note that Razborov's results \cite{Razborov} were crucial for this proof. We also note that we used not only an algorithm for solving systems of equations in a free group, but also the fact (due to \cite{Razborov} as well) that it is algorithmically possible to find a subgroup of $F_n$ of the maximum rank generated by a solution tuple of the given system of equations.

\section{Proof of Theorem \ref{limit}}

%The proof goes essentially along the same lines as the proof of Theorem \ref{free} does.
For the most part, the proof is similar to that of Theorem \ref{free}.
Again, let $g_1,\ldots, g_m$ be the given generators of the group $G$, and Let $r_1,\ldots, r_s$ be all defining relators of $G$. Let $\alpha : g_i \to x_i$ for some $x_i \in H$. This map extends to a homomorphism $\alpha : G \to H$ if and only if $\alpha(r_i)=1$ for all $i=1, \ldots, s$. This translates into a system of $s$ equations in the group $H$.

There are known algorithms for solving systems of equations in limit groups (see e.g. \cite{OMRS}). Moreover, the results of \cite{OMRS} imply that in a limit group $H$, different $m$-tuples of solutions of a system of equations generate only finitely many subgroups $H_i$ of the group $H$ up to isomorphism, and a (finite) presentation
of each subgroup $H_i$ can be algorithmically computed according to \cite[Theorem 30]{OKAM2006}.

We will therefore run an algorithm from \cite{OMRS} to solve the system of equations mentioned in the first paragraph of this section. Then we
find generating $m$-tuples $(h_{i1},\ldots, h_{im})$ of subgroups $H_i$. Then, using an algorithm from \cite{OMRS}, we find (finitely many) defining relations for each subgroup $H_i$ representing an isomorphism class mentioned in the previous paragraph.

Thus, if $G$ can be embedded in $H$, it should be isomorphic to one of the subgroups $H_i$. Suppose there are $k$ of them. We will then run $k$ algorithms $\mathcal{C}_i$ in parallel, where each $\mathcal{C}_i$, in turn, is a pair of algorithms $(\mathcal{A}_i, \mathcal{B}_i)$ running in parallel.

As in the proof of Proposition \ref{main}, algorithm $\mathcal{A}_i$ will detect non-isomorphism by looking for a nontrivial element in the kernel of $\varphi: g_j \to h_{ij}$.
If the kernel is nontrivial, then the subgroup $H_i$ is isomorphic to a proper factor group of the group $G$ and therefore cannot be isomorphic to $G$ itself because all finitely generated subgroups of a limit group are Hopfian.

At the same time, algorithm $\mathcal{B}_i$ will detect isomorphism of the subgroup $H_i$ to the group $G$ by looking for a map $\psi$, given on the generators $h_{ij}$ of $H_i$, such that $\psi(\varphi(g_i))=g_i$ for all generators $g_i$ of the group $G$. This is done the same way as in the proof of Proposition \ref{main}, but there is one more ingredient needed here. To check if $\psi$ is a homomorphism, we see if $\psi$ takes each defining relation of $H_i$ to the identity element of $G$.

Eventually one of the algorithms, $\mathcal{A}_i$ or $\mathcal{B}_i$, will stop and give an answer about isomorphism (or non-isomorphism) of $H_i$ to $G$.

\bigskip

\noindent {\bf Acknowledgement.} I am grateful to Olga Kharlampovich and Alexei Myasnikov for useful discussions on equations in groups and on various properties of limit groups.

\end{document}